\documentclass[11pt,a4paper]{article}
\ifdefined\pdfminorversion\pdfminorversion=7\fi

\usepackage{fancyhdr}
\usepackage{amsmath}
\usepackage{bm}
\usepackage{esint}
\usepackage{amsfonts}
\usepackage{amsthm}
\usepackage{amssymb}
\usepackage{amsbsy}
\usepackage{verbatim}

\usepackage{url}
\usepackage{hyperref}
\usepackage{mathrsfs}
\usepackage[hmargin = 1in,vmargin=.75in]{geometry}
\usepackage{color}
\usepackage{amstext}
\usepackage{enumerate}
\usepackage{pifont}
\usepackage{prettyref}
\usepackage{todonotes}
\usepackage{xcolor}

\usepackage{graphicx}
\usepackage{subcaption} 
\usepackage{epstopdf}   
\usepackage{booktabs}
\usepackage{needspace}
\usepackage{placeins}

\font\msbm=msbm10

\numberwithin{equation}{section}

\theoremstyle{plain}
\newtheorem{theorem}{Theorem}[section]
\newtheorem{lemma}[theorem]{Lemma}
\newtheorem{corollary}[theorem]{Corollary}

\newtheorem{proposition}[theorem]{Proposition}
\newtheorem{definition}{Definition}[section]

\def\mathbb#1{\hbox{\msbm{#1}}}

\newcommand{\bt}{\boldsymbol{t}}

\newcommand{\bv}{\boldsymbol{v}}

\newcommand{\BA}{\boldsymbol{A}}

\newcommand{\BD}{\boldsymbol{D}}
\newcommand{\BE}{\boldsymbol{E}}

\newcommand{\I}{\boldsymbol{I}}

\newcommand{\BO}{\boldsymbol{O}}

\newcommand{\BQ}{\boldsymbol{Q}}
\newcommand{\BR}{\boldsymbol{R}}

\newcommand{\BU}{\boldsymbol{U}}
\newcommand{\BV}{\boldsymbol{V}}

\newcommand{\BX}{\boldsymbol{X}}
\newcommand{\BY}{\boldsymbol{Y}}
\newcommand{\BZ}{\boldsymbol{Z}}

\newcommand{\BDelta}{\boldsymbol{\Delta}}
\newcommand{\BPhi}{\boldsymbol{\Phi}}

\newcommand{\BSigma}{\boldsymbol{\Sigma}}

\newcommand{\BLambda}{\boldsymbol{\Lambda}}

\newcommand{\bzero}{\boldsymbol{0}}

\newcommand{\RR}{\mathbb{R}}
\newcommand{\PP}{\mathcal{P}}

\newcommand{\lag}{\langle}
\newcommand{\rag}{\rangle}

\DeclareMathOperator{\Tr}{Tr}
\newcommand{\eps}{\epsilon}

\newcommand*\diff{\mathop{} \mathrm{d}}

\DeclareMathOperator{\St}{St}

\DeclareMathOperator{\Od}{O}

\DeclareMathOperator{\E}{\mathbb{E}}

\DeclareMathOperator{\diag}{diag}

\DeclareMathOperator{\Ran}{Ran}

\DeclareMathOperator{\blkdiag}{blkdiag}
\DeclareMathOperator{\rank}{rank}
\DeclareMathOperator{\argmin}{argmin}
\DeclareMathOperator{\argmax}{argmax}

\renewcommand{\Pr}{\mathbb{P}}

\DeclareMathOperator{\SO}{SO}

\usepackage{array,tabularx}
\hypersetup{hidelinks}
\newif\ifshowrevisions
\showrevisionstrue
\ifshowrevisions
  \definecolor{revisionblue}{RGB}{0,55,165}
\else
  \definecolor{revisionblue}{RGB}{0,0,0}
\fi

\begin{document}

\title{\bf Statistical Analysis of the Structure from Noise Phenomenon in Generalized Procrustes Problems}

\date{\today}

\author{Jiayang Yin\thanks{Center for Data Science, New York University Shanghai, China. S.L. and J.Y. are (partially) financially supported by the National Key R\&D Program of China, Project Number 2021YFA1002800, National Natural Science Foundation of China No.12571105, Shanghai STCSM Rising Star Program No. 24QA2706200, STCSM General Program No. 24ZR1455300, SMEC AI Initiative Program, and NYU Shanghai Boost Fund. J.Y. is also supported by NYU Shanghai Ph.D. fellowship.}~\thanks{Center for Data Science, New York University.},~ Shuyang Ling$^*$}

\maketitle

\begin{abstract}
In cryo-electron microscopy (cryo-EM), observations are highly noisy, and scientists often use reference templates to guide matching, alignment, and averaging to reconstruct the molecular structure. However, these procedures can produce an output resembling the template even when the observations contain only noise, a phenomenon known as Einstein from noise (EfN). This phenomenon has attracted attention in statistics, where simplified observation models involving group actions have been studied to understand the resulting model bias. Characterizing the bias in the full cryo-EM pipeline remains challenging. We study template matching and averaging in the generalized orthogonal Procrustes problem, which provides an analytically tractable setting involving continuous orthogonal group action and lower-dimensional projections of ordered point clouds. We show that, for pure Gaussian noise, the mean of the aligned-and-averaged estimator admits the same left and right singular-vector bases as the template but its singular values may be distorted. This identifies the subspace structure imposed by the template even in the absence of signal. We characterize the singular values and bound their distortion, establish finite-sample operator-norm bounds, and give sufficient sample sizes for the population mean to dominate random fluctuations. In fact, this bias still exists under weak signals: when the template rank is at least the observation dimension, the change from the pure-noise mean is $O(\eta^2)$ in signal amplitude $\eta$. Our results provide an explicit characterization of model bias in a registration problem motivated by cryo-EM and show why averaging more observations can stabilize structure imposed by the template without removing the bias.
\end{abstract}

\noindent{\bf Keywords:} structure from noise; generalized Procrustes problem; template matching; model bias; cryo-EM.

\noindent{\bf MSC2020:} 62H35, 62H12, 60B20.

\section{Introduction}

In single-particle cryo-electron microscopy (cryo-EM), one observes highly noisy two-dimensional projections of a molecular structure at unknown locations and orientations~\cite{BBS20}. The aim is to recover the three-dimensional structure from these observations. The electron exposure must be limited to reduce radiation damage, so the structural information in an individual image can be weak. Therefore, scientists combine many observations and often use reference templates to guide particle selection, matching, and alignment. The images assigned to a common view can then be averaged to enhance the signal and the information from different views is combined in three-dimensional reconstruction. The success of this approach depends on whether the processing extracts information from the observations or introduces structure through the chosen reference.

The \emph{Einstein-from-noise} (EfN) experiment illustrates this concern in a simple setting~\cite{SHBG09,H13}. Suppose that the observations are believed to be noisy, shifted copies of a reference image, such as a portrait of Einstein, but in fact contain only independent noise. Each observation is aligned to the reference by maximizing their cross-correlation over shifts, and the aligned observations are averaged. Remarkably, the output resembles the portrait even though it is absent from every observation. Averaging the original centered noise would converge to zero. However, after alignment,  the transformation applied to each observation depends on that same noise, and the transformed observations can have a nonzero population mean. Thus increasing the number of observations can make the reference-like structure more apparent.

This impact has long been recognized in structural biology. Stewart and Grigorieff~\cite{SG04} studied noise bias in single-particle refinement. Henderson~\cite{H13} explained how template-based particle selection can produce apparently meaningful three-dimensional maps from noise, and van Heel~\cite{vH13} demonstrated the effect using random-noise micrographs. These examples show why the output of a reconstruction procedure requires careful interpretation, especially when the signal is too weak to distinguish individual particles reliably. We use \emph{structure from noise} to describe the broader appearance of structure in noise-only data, with EfN as a canonical example. The related term \emph{confirmation bias} refers to the influence of prior expectations, encoded through templates or other modeling choices, on the inferred structure.

These observations motivate a statistical question: which features of an estimated structure are inherited from the reference, and which are supported by the data? A quantitative answer for the full cryo-EM pipeline is difficult because particle selection, pose estimation, classification, and reconstruction interact. Simplified observation models involving group actions provide a way to isolate the effect of alignment. By specifying how an unknown transformation acts on a signal, one can study how fitting that transformation to noise changes the output of the estimation procedure.

Balanov, Huleihel, and Bendory~\cite{BHB24} analyze this mechanism for fixed-template alignment over cyclic shifts. For white Gaussian noise and under their template assumptions, they show that the Fourier phases of the aligned average converge to those of the reference and quantify the convergence. This identifies a structural feature that survives the alignment-and-averaging procedure and helps explain why the output resembles the template even when their Fourier magnitudes differ. Their analysis suggests a natural question for other transformation models: what plays the role of the preserved Fourier phases, and how does the transformation model determine the resulting bias?

A complementary line of work studies the bias introduced before alignment, when observations are selected using a reference. Balanov, Zabatani, and Bendory~\cite{BZB25} analyze template-based particle picking and its effects on downstream class-mean estimation and volume reconstruction under specified models. Their selection mechanism differs from aligning a given set of observations, although both can produce structure from noise. This distinction is relevant to the present work where every observation is matched to a prescribed template and included in the average. Our focus is the bias created by fitting transformations to the observations themselves.

\paragraph{The Procrustes model.}
In this work, we study template matching followed by averaging in the generalized orthogonal Procrustes problem. The observations are ordered point clouds related to an unknown signal by orthogonal transformations and projection to a lower-dimensional space. For a prescribed template, each pose is estimated by maximizing its correlation with the observation. The observations are then mapped back to the template coordinates and averaged. This model incorporates continuous transformations, including rotations, and projected observations, which are relevant to cryo-EM. At the same time, known point correspondences and the explicit SVD solution of the matching step make the dependence between the fitted pose and the noise accessible to analysis. 

Our main question is: \emph{which features of the prescribed template are preserved in the mean of the structure-from-noise estimator when the observations contain only noise?} We also ask how large the resulting bias is, how many observations are needed for it to dominate sampling fluctuations, and whether it persists when a weak signal is present. These questions connect the qualitative resemblance observed in EfN experiments to quantities that can be analyzed in the Procrustes model.

\paragraph{Main results and implications.}
If the observation is purely independent standard Gaussian noise, we show that the expected structure-from-noise estimator admits the same left and right singular-vector bases as the template, while its singular values may change. Thus the template determines the directional structure of the population output even in the absence of signal. Geometrically, the mean reproduces the template after rescaling along its principal directions. This is analogous to Fourier-phase preservation in cyclic-shift EfN, with singular-vector structure providing the preserved feature in the Procrustes setting. The statement concerns common singular-vector bases and allows for their nonuniqueness when singular values are repeated. Our main results are as follows.

\begin{enumerate}
\item {\bf Population structure and spectral distortion.} We characterize the mean of the structure-from-noise estimator and bound its singular values in terms of the template spectrum, dimensions, and rank. When the nonzero template singular values are equal, the mean is a positive scalar multiple of the template, with an explicit expression for the scalar (Theorem~\ref{thm:main} and Corollary~\ref{cor:svd_efn}).

\item {\bf Finite-sample behavior.} We establish operator-norm bounds for deviations from the  mean and give sufficient sample sizes for this mean to dominate random fluctuations. The dependence on the dimensions and template rank makes the required sample size explicit (Theorem~\ref{thm:error}).

\item {\bf Weak-signal scenario.} We expand the population mean of the structure-from-noise estimator around zero signal strength. When the template rank is at least the observation dimension, the linear term vanishes and the change from the pure-noise mean is $O(\eta^2)$ in signal amplitude $\eta$. Thus a weak signal can still preserve the leading structure of the pure-noise bias (Theorem~\ref{thm:main_lsnr}).
\end{enumerate}

Together, these results identify both the structure of the bias and the regimes in which it dominates the observed output. Increasing the sample size reduces fluctuations around the population mean but does not remove its dependence on the template. A stable average resembling the reference can  arise without an underlying signal, and this effect persists under weak signals. By making the preserved directional structure and its spectral distortion explicit, the Procrustes model provides a tractable analysis of an alignment mechanism that may help explain model bias in more complicated cryo-EM procedures.

\paragraph{Related statistical and geometric analyses.}
Other statistical models expose complementary aspects of the same concern. Wang et al.~\cite{WYCT21} study spherical noise vectors selected for high correlation with a reference while Gaussian-mixture analyses examine bias in class-mean estimation~\cite{BBH25}. Balanov, Huleihel, and Bendory~\cite{BHB26} analyze initialization bias and finite-sample EfN behavior for low-SNR expectation-maximization; their iterative soft assignments differ from fixed-template hard alignment. Xu et al.~\cite{XBSB26} study Bayesian orientation estimation over $\SO(3)$ and report empirical mitigation of EfN artifacts through posterior averaging, illustrating the importance of the pose-estimation rule.

Procrustes registration provides a tractable setting for studying continuous transformations. Its classical orthogonal least-squares solution~\cite{S66} underlies a statistical framework for landmark shape analysis~\cite{G91}. Consistency for an underlying shape depends on the geometry and noise assumptions: Kent and Mardia~\cite{KM97} establish planar shape-consistency results under isotropic errors while Huckemann~\cite{H11} distinguishes convergence to a population Procrustes mean from agreement with the shape of an underlying perturbation model. More generally, Devilliers and coauthors~\cite{DATP17a,DATP17b} establish inconsistency and quantify bias in unknown-template estimation under group actions, including zero-signal cases. Miolane, Holmes, and Pennec~\cite{MHP17} derive a geometric small-noise bias expansion. 
These results establish that alignment-related population bias can persist with arbitrarily many observations. While these works study the consistency of an estimated template, we characterize which features of a prescribed reference are preserved or distorted in the aligned average of pure noise.

A closely related observation model is studied by Pumir, Singer, and Boumal~\cite{PSB21}, who recover an unknown point cloud from noisy orthogonal copies using Gram-matrix invariants when the poses cannot be reliably estimated. Recovery guarantees for generalized Procrustes estimation and multi-reference alignment also appear in~\cite{CKS15,L23b,L25,PWB+19,RBO21}. We instead characterize the bias of a fixed-template alignment procedure. This complements signal-recovery theory by identifying the structure that this procedure can produce in pure noise and retain in the presence of a weak signal.
\subsection{Organization}
Section~\ref{sec: preliminaries} defines the observation model and estimator and explains their connection to the cryo-EM pipeline. Section~\ref{s:main} presents the main results and numerical experiments, and Section~\ref{s:proof} contains the proofs.

\subsection{Notation}\label{ss:notation}
Bold uppercase and lowercase letters denote matrices and vectors; $\I_d$ is the $d\times d$ identity. We use $\|\cdot\|$, $\|\cdot\|_*$, and $\|\cdot\|_F$ for the operator, nuclear, and Frobenius norms, and $\langle\BX,\BY\rangle=\Tr(\BX^{\top}\BY)$. The symbol  $\BX^\dagger$ denotes the Moore--Penrose pseudoinverse. We use the positive semidefinite order $\preceq$, positive semidefinite square roots, and standard asymptotic notation. Finally, $a\vee b=\max\{a,b\}$, $a\wedge b=\min\{a,b\}$, and $[n]=\{1,\cdots,n\}$.

\section{Preliminaries}\label{sec: preliminaries}

\subsection{Observation model and estimators}\label{sec: generalized_procrustes}
Let $\BA^{\star}\in\RR^{d\times m}$ be a fixed unknown signal whose columns represent $m$ ordered points in $\RR^d$. We observe $n$ noisy, possibly projected copies of this point cloud:
\[
\BA_i=\BR_i^{\star}\BA^{\star} + \BDelta_i\in\RR^{r\times m}, \qquad i\in[n],\qquad 1\leq r\leq d\leq m,
\]
where the additive noise $\{\BDelta_i\}_{i=1}^n$ are independent matrices with i.i.d. standard Gaussian entries. The hidden poses have orthonormal rows:
\[
\BR_i^{\star}\in\St(r,d) :=\{\BR\in\RR^{r\times d}:\BR\BR^{\top}=\I_r\}.
\]
In particular, we have $\St(d,d)=\Od(d)$, and $\SO(d)=\{\BR\in\Od(d):\det\BR=1\}$.

The poses may be deterministic or independent of the noise. When $r=d$, they are orthogonal matrices. When $r<d$, each $\BR_i^{\star}$ is a flat rectangular matrix representing an orthogonal transformation followed by projection onto $r$ coordinates. In particular, $d=3$ and $r=2$ model the loss of one spatial dimension, providing a simple analogue of the lower-dimensional observations in cryo-EM. Thus $d$, $r$, $m$ and $n$ denote the ambient dimension, the observation dimension,  the number of points, and the number of observations respectively.

The orthogonal Procrustes problem aligns two matrices, or two point clouds with known correspondences, by a rigid orthogonal transformation~\cite{S66}. Here, given a fixed prescribed template $\BA^{\#}\in\RR^{d\times m}$, we match each observation to the template by maximizing their correlation:
\begin{equation}\label{eq: original_procrustes_model}
\widehat{\BR}_i\in\argmax_{\BR\in\St(r,d)} \langle\BA_i,\BR\BA^{\#}\rangle =\argmax_{\BR\in\St(r,d)} \langle\BR,\BA_i\BA^{\#\top}\rangle.
\end{equation}

When $r=d$, this is the usual orthogonal Procrustes problem, equivalent to minimizing $\|\BA_i-\BR\BA^{\#}\|_F^2$. When $r<d$, the term $\|\BR\BA^{\#}\|_F^2$ varies with $\BR$, so the correlation-matching rule is generally different from least squares or Gaussian maximum likelihood. The model assumes known point correspondences and does not estimate permutations or translations. To solve~\eqref{eq: original_procrustes_model}, we introduce the polar retraction as follows.

\begin{definition}[Polar retraction]\label{def: polar_retraction}
For $\BX\in\RR^{r\times d}$ with $r\leq d$, a polar retraction is a maximizer
\[
\PP(\BX)\in\argmax_{\BR\in\St(r,d)}\langle\BR,\BX\rangle
=\argmin_{\BR\in\St(r,d)}\|\BR-\BX\|_F^2.
\]
We use the same definition for tall matrices, with orthonormal columns. A thin SVD $\BX=\BU_{\BX}\BSigma_{\BX}\BV_{\BX}^{\top}$, including zero singular directions, gives a solution $\PP(\BX)=\BU_{\BX}\BV_{\BX}^{\top}$. It is unique at full row or column rank; otherwise the null-space completion can vary.
\end{definition}

Throughout, we use a fixed measurable maximizing rule depending only on its matrix argument, with paired choices under transposition, $\PP(\BX^{\top})=\PP(\BX)^{\top}$. The same rule is used for every observation and every signal level.

\begin{proposition}[Polar identities]\label{prop: full_rank_case_polar_retraction}\label{prop:polar}
If $\BX$ is full row rank, $\PP(\BX)=(\BX\BX^{\top})^{-1/2}\BX$; if $\BX$ is full column rank, $\PP(\BX)=\BX(\BX^{\top}\BX)^{-1/2}$. Every polar retraction satisfies
\[
\PP(\BX)\BX^{\top}=(\BX\BX^{\top})^{1/2}, \qquad \PP(\BX)^{\top}\BX=(\BX^{\top}\BX)^{1/2}.
\]

\end{proposition}

The polar identities follow directly from the SVD. The fitted pose in~\eqref{eq: original_procrustes_model} is $\widehat{\BR}_i=\PP(\BA_i\BA^{\#\top})$. We define the alignment and backprojection map by
\[
{\cal T}(\BY):=\PP(\BY\BA^{\#\top})^{\top}\BY\in\RR^{d\times m}, \qquad \BY\in\RR^{r\times m}.
\]
For $r=d$, the transpose of the fitted pose reverses an orthogonal transformation. For $r<d$, it maps the observation back to the template coordinates, but $\widehat{\BR}_i^{\top}\widehat{\BR}_i$ is only a rank-$r$ projection and cannot restore the missing components. The aligned-and-averaged estimator is
\[
\widehat{\BA}
=\frac1n\sum_{i=1}^n{\cal T}(\BA_i)
=\frac1n\sum_{i=1}^n\widehat{\BR}_i^{\top}\BA_i.
\]

\paragraph{The structure-from-noise estimator.}
Setting $\BA^{\star}=\bzero$ gives purely Gaussian observations $\BA_i=\BDelta_i$. Applying the same matching and averaging procedure yields the structure-from-noise (SfN) estimator
\begin{equation}\label{eq: estimator_EfN}
\widehat{\BA}_{\rm SfN}:=\frac1n\sum_{i=1}^n{\cal T}(\BDelta_i).
\end{equation}
Let $k=\rank(\BA^{\#})$ and if $k\geq 1$, denote its ordered nonzero singular values by $\sigma_1\geq\cdots\geq\sigma_k>0$. For Gaussian noise, $\rank(\BDelta_i\BA^{\#\top})=r\wedge k$ holds almost surely. The fitted pose is unique almost surely when $k\geq r$; otherwise the fixed maximizing rule above specifies the estimator.

If the fitted pose were independent of the noise, its backprojection would have zero mean. Fitting it to the same observation instead gives, 
\[
\begin{aligned}
\left\langle\E\widehat{\BA}_{\rm SfN},\BA^{\#}\right\rangle = \E\max_{\BR\in\St(r,d)}\langle\BDelta_i,\BR\BA^{\#}\rangle =\E\|\BDelta_i\BA^{\#\top}\|_* >0
\end{aligned}
\]
for any $\BA^{\#} \neq 0.$ Thus template fitting creates a nonzero mean from centered noise. Theorem~\ref{thm:main} introduces the template SVD and characterizes the structure of this mean.

\subsection{Connection to the cryo-EM pipeline}

For a particle image extracted from a micrograph, a  cryo-EM forward model is~\cite{BBS20}
\begin{equation}\label{def:cryoem}
f_i\approx C_i\tau_{\bt_i}\Pi_{\BR_i}\mu+\eps_i.
\end{equation}
Here $\mu$ is the molecular density, $\BR_i\in\SO(3)$ the viewing orientation, $\Pi_{\BR_i}$ line integration through the rotated density, $\tau_{\bt_i}$ an in-plane translation, and $C_i$ the contrast transfer function (CTF) operator. These images are noisy lower-dimensional observations whose unknown poses need to be estimated when combining information from different particles.

The Procrustes model retains two ingredients of this setting: lower-dimensional views of an unknown object and template-dependent pose fitting followed by aggregation. Its flat row-orthogonal matrices mimic the loss of a spatial dimension, while $\widehat{\BR}_i^{\top}\BA_i$ is a simple backprojection statistic. In both settings, fitting the pose to the same noisy observation can favor agreement with the reference and produce structure after averaging. Known point correspondences and the explicit SVD solution of the matching step make this dependence tractable in the Procrustes model.

The analogy concerns this alignment-and-averaging mechanism. Projecting corresponding points differs from line integration of a density in cryo-EM, and $\Od(d)$ includes reflections whereas physical orientations belong to $\SO(3)$. Our observations are already given; class assignments, translations, CTFs, and the full tomographic inverse problem are not modeled. The resulting proxy allows us to characterize alignment-induced bias explicitly while keeping it distinct from selection bias introduced during particle picking.

\section{Main results}\label{s:main}

We study three aspects of the structure-from-noise estimator:
\begin{enumerate}
\item Which features of the prescribed template are preserved in the population mean $\E\widehat{\BA}_{\rm SfN}$?
\item How many observations are sufficient for this mean to dominate random fluctuations?
\item How does the mean change when a weak signal is present?
\end{enumerate}
We begin with the singular-vector structure of the pure-noise population mean.

\begin{theorem}[\bf Population structure of the SfN estimator]\label{thm:main}
Fix $1\leq r,k\leq d\leq m$, and let a deterministic template of rank $k$ have the SVD
\[
\BA^{\#}=\BU^{\#}\BSigma^{\#}\BV^{\#\top},
\qquad \BU^{\#}\in\Od(d),\quad \BV^{\#\top}\in\St(d,m),
\]
where $\BSigma^{\#} = \diag(\sigma_1,\ldots,\sigma_k,0,\ldots,0)$ and $\sigma_1\ge\cdots\ge\sigma_k>0$. For independent standard Gaussian observations $\BA_i=\BDelta_i\in\RR^{r\times m}$, the estimator in~\eqref{eq: estimator_EfN} satisfies

\begin{equation}\label{eq: population_mean_EfN}
\E\widehat{\BA}_{\rm SfN} = \BU^{\#}\BSigma_{\rm SfN}\BV^{\#\top}, ~~~\BSigma_{\rm SfN} :=\E\left[\PP(\BPhi\BSigma^{\#})^{\top}\BPhi\right],
\end{equation}
where $\BPhi\in\RR^{r\times d}$ is standard Gaussian and $\BSigma_{\rm SfN}$ is diagonal, with positive first $k$ entries and zero remaining entries. More precisely, it holds that
\[
\BSigma_{\rm SfN}  =  \BLambda\BSigma^{\#\dagger} = \blkdiag\big(\BLambda_+\BSigma^{\#-1}_+, \bzero_{(d-k)\times(d-k)}\big)
\]
where $\BSigma_+^{\#}= \diag(\sigma_1,\cdots,\sigma_k)\in\RR^{k\times k}$, $\BPhi_+$ is the first $k$ columns of $\BPhi,$ $\BLambda = \E\left[\PP(\BPhi\BSigma^{\#})^{\top}\BPhi\BSigma^{\#}\right]$ and $\BLambda_+ = \E(\BSigma_+^{\#}\BPhi_+^{\top}\BPhi_+\BSigma_+^{\#})^{1/2}$.
Therefore, the population mean of $\widehat{\BA}_{\rm SfN}$ and the template admit the same left and right singular-vector bases.
\end{theorem}

Theorem~\ref{thm:main} shows that the population mean reproduces the template after rescaling along its principal directions. The statement concerns common singular-vector bases, allowing for their non-uniqueness when singular values repeat. The next result quantifies the spectral distortion and gives an exact expression when the nonzero template singular values are equal. 

\begin{corollary}[Singular-value bounds for the population mean]\label{cor:svd_efn}
Under the assumptions of Theorem~\ref{thm:main}, let $\BPhi_+\in\RR^{r\times k}$ be standard Gaussian. The population mean $\E\widehat{\BA}_{\rm SfN}\in\RR^{d\times m}$ has the spectral bounds
\begin{equation}\label{eq: singular_value_bounds_Lambda}
\frac{\E\|\BPhi_+\|_*}{k\sigma_1}\BSigma^{\#}
\preceq\BSigma_{\rm SfN}
\preceq\frac{\E\|\BPhi_+\|_*}{k\sigma_k}\BSigma^{\#}.
\end{equation}
Equivalently, its nonzero singular values satisfy
\begin{equation}\label{eq: singular_value_bounds_tau}
\frac{\sigma_j}{k\sigma_1}\E\|\BPhi_+\|_*
\leq(\BSigma_{\rm SfN})_{jj}
\leq\frac{\sigma_j}{k\sigma_k}\E\|\BPhi_+\|_*,
\qquad j\in[k],
\end{equation}
and $(\BSigma_{\rm SfN})_{jj}=0$ for $k+1\leq j\leq d$. If $\BSigma^{\#}=\blkdiag(\I_k,\bzero_{(d-k)\times(d-k)})$, the lower and upper bounds coincide:
\begin{equation}\label{eq: singular_value_isotropic_EfN}
\BSigma_{\rm SfN}=\blkdiag\left(\frac{\E\|\BPhi_+\|_*}{k}\I_k,\bzero_{(d-k)\times(d-k)}\right).
\end{equation}
In general, the operator norm of the mean is bounded below by
\[
\|\E\widehat{\BA}_{\rm SfN}\| \geq\frac{\E\|\BPhi_+\|_*}{k} \geq\frac{r}{\sqrt{r+k+1}}.
\]
\end{corollary}

The Gaussian nuclear norm determines the scale of the population mean. Corollary~\ref{cor:svd_efn} gives the lower bound
\[
\E \|\BPhi_+ \|_* \geq \frac{rk}{\sqrt{r+k+1}}.
\]
We next evaluate this quantity numerically for small dimensions.

\paragraph{Numerical evaluation of $\E\|\BPhi_+\|_*$}
For $k=d$ and $1\leq r\leq d$, let $\BPhi_+\in\RR^{r\times d}$ be standard Gaussian. The Wishart eigenvalue density gives an integral expression for its expected nuclear norm. The ordered eigenvalues of $\BPhi_+\BPhi_+^{\top}$ lie in
\[
\mathcal{D}_r:=\{(\lambda_1,\cdots,\lambda_r): \lambda_1>\cdots>\lambda_r>0\}.
\]
Their joint density on $\mathcal{D}_r$ is
\[
p(\lambda_1,\cdots,\lambda_r) = K_{r,d}\prod_{1\leq i<j\leq r}(\lambda_i-\lambda_j) \prod_{i=1}^{r}\lambda_i^{(d-r-1)/2}e^{-\lambda_i/2},
\]
where
\[
K_{r,d} =\frac{\pi^{r/2}}{2^{rd/2} \displaystyle\prod_{i=1}^{r} \Gamma(i/2)\Gamma((d-i+1)/2)}.
\]
This is the real Wishart density in~\cite[Section~4.5, Eqs.~(4.5)--(4.6)]{ER05}, with the normalization multiplied by $r!$ to account for ordering; see also~\cite{M91}. Since the singular values of $\BPhi_+$ are $\sqrt{\lambda_1},\cdots,\sqrt{\lambda_r}$, we obtain
\begin{equation}\label{eq: Gaussian_constant_integral}
c_{r,d}:=\frac{1}{d} \E \|\BPhi_+\|_* =\frac{1}{d}\int_{\mathcal{D}_r}
\left(\sum_{i=1}^{r}\sqrt{\lambda_i}\right) p(\lambda_1,\cdots,\lambda_r) {\rm d}\lambda_1\cdots {\rm d}\lambda_r.
\end{equation}
The ordered integration domain is essential for this normalization of $p$. Table~\ref{tab:nuclear_norm_expectation} reports numerical values for dimension choices motivated by applications in computer vision, structural biology, and robotics.

\begin{table}[h!]
\centering
\caption{Numerical values of $c_{r,d}=d^{-1} \E \|\BPhi_+\|_*$ for $\BPhi_+\in\mathbb{R}^{r\times d}$ with independent standard Gaussian entries, rounded to six decimal places.}
\begin{tabular}{c|c|l}
\toprule
$(r,d)$ & $c_{r,d}$ & Illustrative application \\
\midrule
$(2,2)$ & $1.145797$ & 2D point-cloud registration (computer vision) \\
$(2,3)$ & $1.008589$ & Projected point-cloud model (3D to 2D) \\
$(3,3)$ & $1.418150$ & 3D point-cloud registration (robotics and computer vision) \\
\bottomrule
\end{tabular}
\label{tab:nuclear_norm_expectation}
\end{table}

We now quantify fluctuations around the population mean and give sufficient sample sizes for this mean to dominate the random error.

\begin{theorem}[Concentration and sample size]
\label{thm:error}
Under the assumptions of Theorem~\ref{thm:main}, suppose $n\ge d\vee2$, $\gamma\ge1$, and $(\gamma+1)\log n\le\sqrt{nd}$. Then
with probability at least $1-O(n^{-\gamma})$,
\begin{equation}\label{eq: order_EfN_estimator}
\left\|\widehat{\BA}_{\rm SfN} - \E \widehat{\BA}_{\rm SfN}\right\| \leq \frac{\sqrt{r} + \sqrt{m-d} + C(\sqrt r+\sqrt d)\sqrt{2\gamma \log n}}{\sqrt{n}}
\end{equation}
where $C>1$ is an absolute constant. In particular, if also 
\[
n > \frac{(\sqrt{r} + \sqrt{m-d} + C(\sqrt r+\sqrt d)\sqrt{2\gamma \log n})^2 (r + k+1)}{r^2}
\]
then $\left\|\widehat{\BA}_{\rm SfN}-\E \widehat{\BA}_{\rm SfN}\right\| < \left\|\E \widehat{\BA}_{\rm SfN}\right\|.$
In particular, if $d = r = k$, then it suffices to have
\[
n > (2 + d^{-1})\left(1 + \sqrt{\frac{m-d}{d}} + 2C\sqrt{2\gamma \log n}\right)^2.
\]
\end{theorem}

Theorem~\ref{thm:error} compares the magnitude of the population mean with that of the random fluctuations. For fixed $r,d,k,\gamma$, the bound gives a sufficient sample size proportional to $m$ when the $m$-dependent term dominates the logarithmic terms. Its sample-size conditions are sufficient conditions for operator-norm dominance; they do not guarantee that every template feature is visually apparent. The following experiments illustrate the population structure and the dependence on sample size.

\paragraph{Numerical simulations}

\begin{figure}[h!]
\centering
\includegraphics[width=\textwidth]{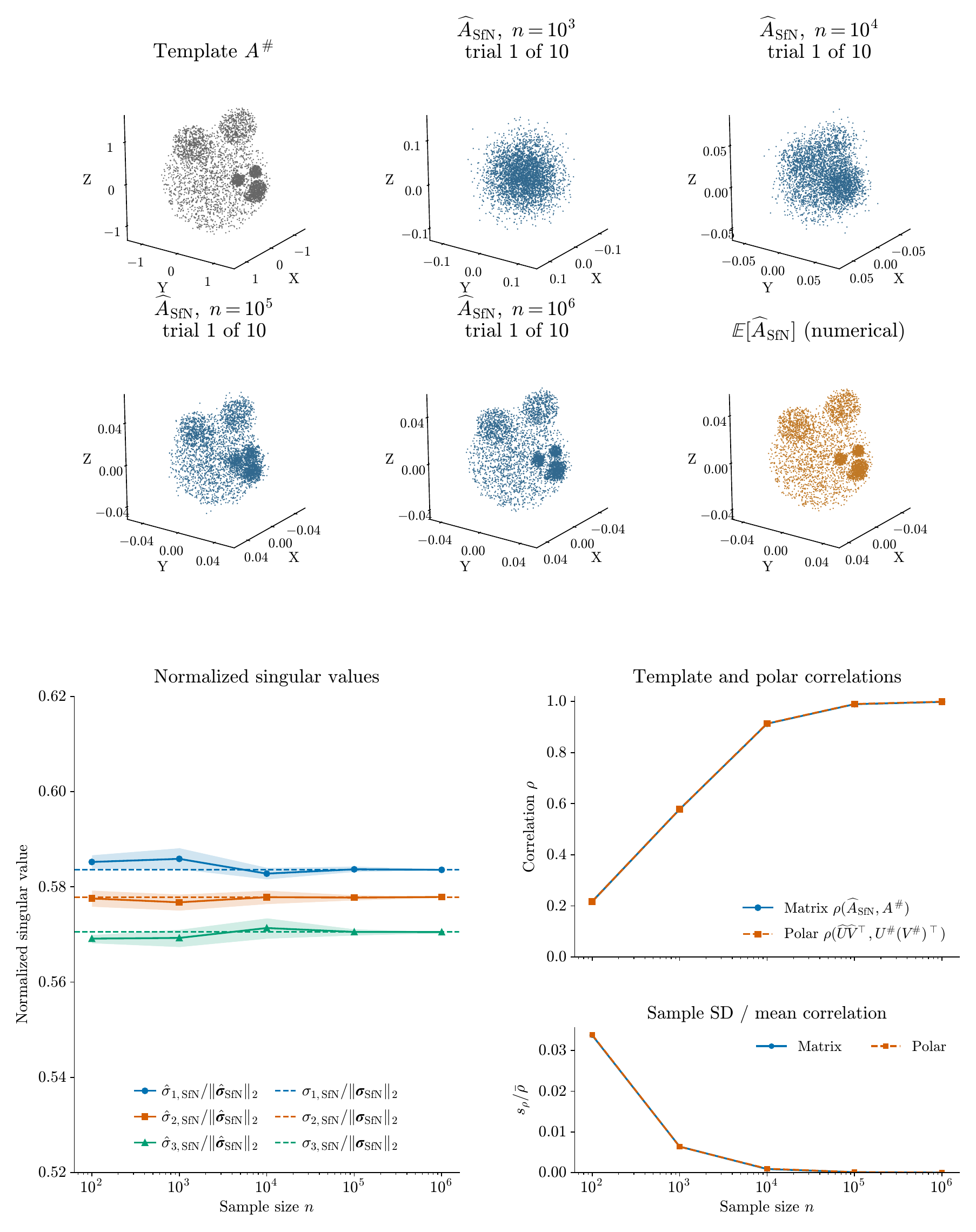}
\caption{SfN estimates without projection as $n$ increases, for $r=d=3$ and $m=4000$. The Mickey mouse head in the $\widehat{\BA}_{\rm SfN}$ becomes visible after averaging over $10^4$ samples, which corresponds to $\rho \approx 0.9.$}
\label{fig:efn_r3}
\end{figure}

\begin{figure}[h!]
\centering
\includegraphics[width=\textwidth]{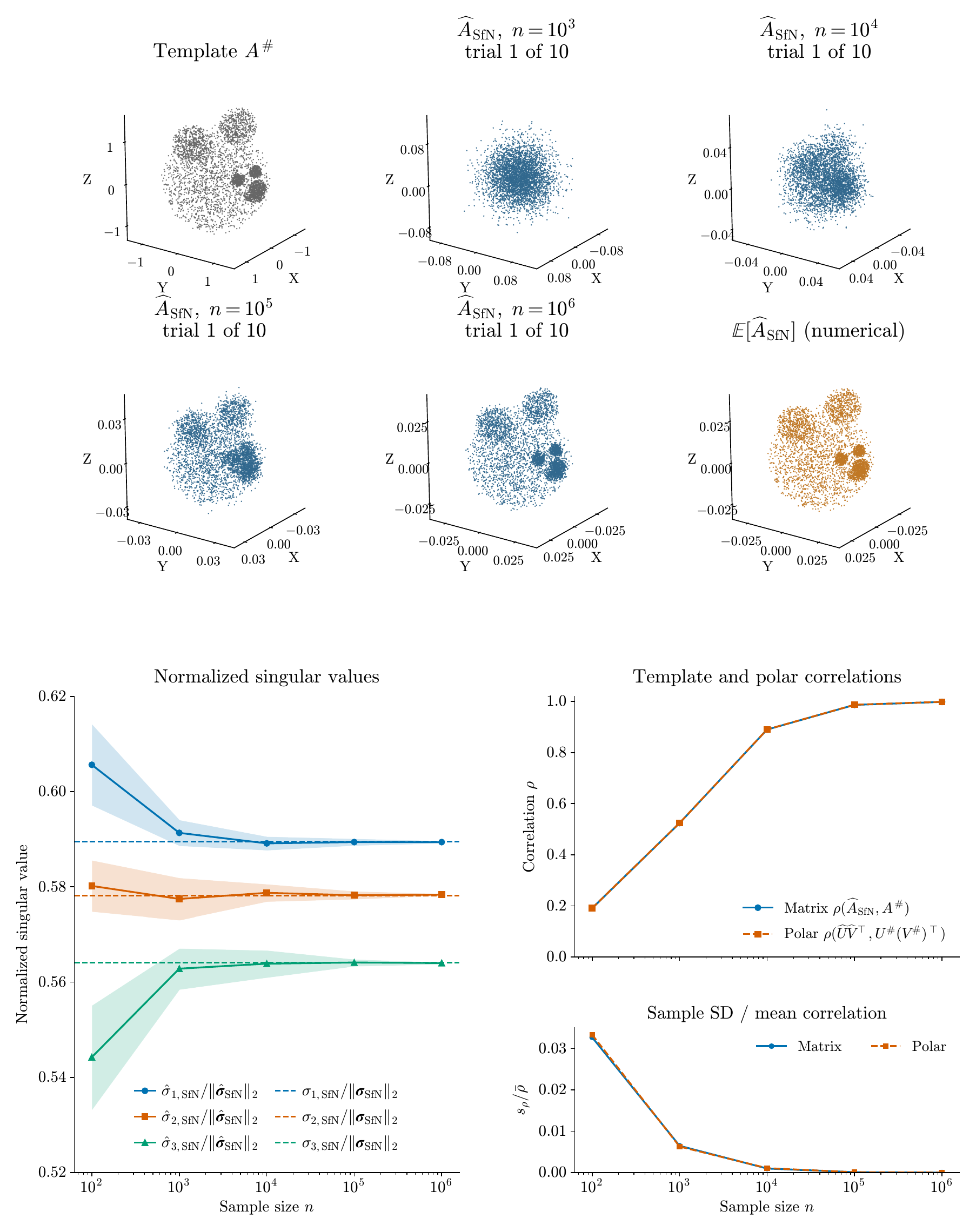}
\caption{SfN estimates with projection as $n$ increases, for $r=2$, $d=3$, and $m=4000$. Similarly as the case of $d= r= 3$, the Mickey mouse head in the $\widehat{\BA}_{\rm SfN}$ becomes visible after averaging over $10^4$ samples, which corresponds to $\rho \approx 0.9.$}
\label{fig:efn_r2}
\end{figure}

Figures~\ref{fig:efn_r3} and~\ref{fig:efn_r2} illustrate the bias for a three-dimensional Mickey Mouse point cloud with $m=4000$ points. Every observation is pure Gaussian noise. We compare full-dimensional observations ($r=d=3$) with projected observations ($r=2$, $d=3$), the latter missing one spatial dimension. In both cases, the estimated shape becomes clearer as $n$ increases, particularly from $n=10^4$ in the displayed examples.

The first six panels show the template, four finite-sample estimates, and the numerically evaluated population mean (second row, right). The bottom-left spectral panel displays the singular values normalized by their Euclidean norm, equivalently by the matrix Frobenius norm: $\widehat\sigma_{j,\rm SfN}/\|\widehat{\BA}_{\rm SfN}\|_F$ and $\sigma_{j,\rm SfN}/\|\E\widehat{\BA}_{\rm SfN}\|_F$, respectively, for $1\le j\le3$. Solid curves average the estimated normalized values over ten independent trials; shaded bands are computed from these trials. The two right panels summarize the alignment scores defined below and their relative variability.

\begin{figure}[h!]
\centering
\includegraphics[width=120mm]{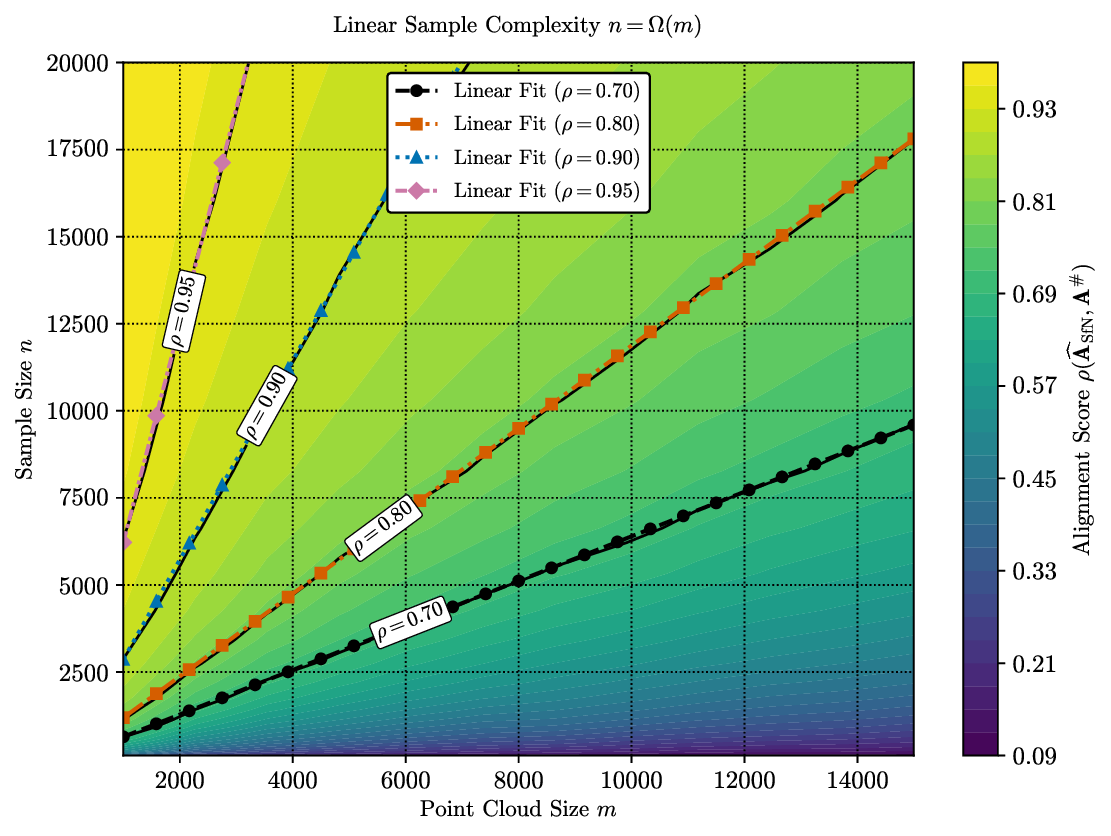}
\caption{Numerical illustration of the sample-size scaling for $d=3$ and $r=2$ with template rank $k\le3$. The alignment-score contours at $\rho\in\{0.70,0.80,0.90,0.95\}$ are compared with linear fits $n=C_{\rho}m$.}
\label{fig:linearitysamplecomplexitycontour}
\end{figure}

To study the joint dependence on the number of points $m$ and observations $n$, Figure~\ref{fig:linearitysamplecomplexitycontour} uses the largest cosine similarity attainable by orthogonally aligning the estimate with the template:
\[
\rho(\widehat{\BA}_{\rm SfN},\BA^{\#}) := \frac{\|\widehat{\mathbf A}_{\rm SfN}\BA^{\#\top}\|_*}{\|\widehat{\mathbf A}_{\rm SfN}\|_F\|\BA^{\#}\|_F} 
\]
Figures~\ref{fig:efn_r3} and~\ref{fig:efn_r2} also provide the polar-factor score $\rho(\PP(\widehat{\BA}_{\rm SfN}),\PP(\BA^{\#}))$. Their variability panels w.r.t. $n$ show the empirical standard deviation divided by the empirical mean of each score over ten trials.
We vary $m$ from $1000$ to $15000$ and $n$ from $100$ to $20000$. The displayed contours are approximately linear and closely follow fits of the form $n=C_{\rho}m$. This supports linear sample-size scaling at the displayed alignment levels over the tested range, consistent with Theorem~\ref{thm:error}.

\paragraph{Persistence under weak signals}
We next ask whether the pure-noise bias persists when a weak signal is present. For the fixed template $\BA^{\#}$ and signal $\BA^{\star}\in\RR^{d\times m}$, consider
\begin{equation}\label{eq:model_lsnr}
\BA_i(\eta)=\eta\BR_i^{\star}\BA^{\star}+\BDelta_i\in\RR^{r\times m},
\quad i\in[n],
\end{equation}
where $\BR_i^{\star}\in\St(r,d)$ is unknown. The parameter $\eta\ge 0$ controls the signal strength. 

We use the same fixed measurable polar-factor selection rule for every observation and every $\eta$.
Then the backprojection estimator is given by
\begin{equation}\label{eq: estimator_LSNR}
\widehat{\BA}(\eta):=\frac{1}{n}\sum_{i=1}^n {\cal T}(\BA_i(\eta)).
\end{equation}
For the same noise realizations, $\widehat{\BA}(0)=\widehat{\BA}_{\rm SfN}$. In particular, $\E\widehat{\BA}(0)=\E\widehat{\BA}_{\rm SfN}=\E{\cal T}(\BDelta_i)$.
The following theorem expands the population mean near $\eta=0$. Expectations also include the poses when they are random. The dimensions, template, signal, and pose distributions are fixed as $\eta\to 0$.

\begin{theorem}[Population expansion in the low-SNR regime]\label{thm:main_lsnr}
Under the model~\eqref{eq:model_lsnr}, for $\eta$ sufficiently small, it holds that
\[
\begin{aligned}
\E[\widehat{\BA}(\eta)] & = \E[\widehat{\BA}_{\rm SfN}]
+\eta\frac{\diff}{\diff\eta}\E\widehat{\BA}(\eta)\Big|_{\eta=0} +\frac{\eta^2}{2}\frac{\diff^2}{\diff\eta^2}\E\widehat{\BA}(\eta)\Big|_{\eta=0}
+O(\eta^3),
\end{aligned}
\]
where
\[
\begin{aligned}
\frac{\diff}{\diff\eta}\E\widehat{\BA}(\eta)\Big|_{\eta=0}
&=\frac{1}{n}\sum_{i=1}^n\E\left[
{\cal T}(\BDelta_i)\langle\BDelta_i,\BR_i^{\star}\BA^{\star}\rangle\right],\\
\frac{\diff^2}{\diff\eta^2}\E\widehat{\BA}(\eta)\Big|_{\eta=0}
&=\frac{1}{n}\sum_{i=1}^n\E\left[
{\cal T}(\BDelta_i)
\bigl(\langle\BDelta_i,\BR_i^{\star}\BA^{\star}\rangle^2
-\|\BR_i^{\star}\BA^{\star}\|_F^2\bigr)\right].
\end{aligned}
\]
If $k\geq r$, the population mean is even in $\eta$, so its first and third derivatives at zero vanish. The expansion reduces to
\begin{equation}\label{eq: LSNR_quadratic_expansion}
\E[\widehat{\BA}(\eta)]
=\E[\widehat{\BA}_{\rm SfN}]
+\frac{\eta^2}{2}\frac{\diff^2}{\diff\eta^2}\E\widehat{\BA}(\eta)\Big|_{\eta=0}
+O(\eta^4).
\end{equation}
In particular, it holds that
\[
\left\|\E[\widehat{\BA}(\eta)]
-\E[\widehat{\BA}_{\rm SfN}]\right\|
=O(\eta^2)
\]
where the remainder bounds are in the operator norm.
\end{theorem}

Theorem~\ref{thm:main_lsnr} concerns the population mean. It provides a perturbation of the estimator $\widehat{\BA}$ in the low-SNR regime. In particular,~\eqref{eq: LSNR_quadratic_expansion} follows from sign symmetry: when $k\geq r$, the matching matrix has full row rank almost surely and the polar factor changes sign with its argument, so the backprojected observation is even in the noise. Expanding the shifted Gaussian density justifies the population expansion without differentiating the polar factor near singular matching matrices. When $k<r$, the general expansion still applies, but cancellation of the linear term is not ensured.

\begin{figure}[h!]
\centering
\includegraphics[width=110mm]{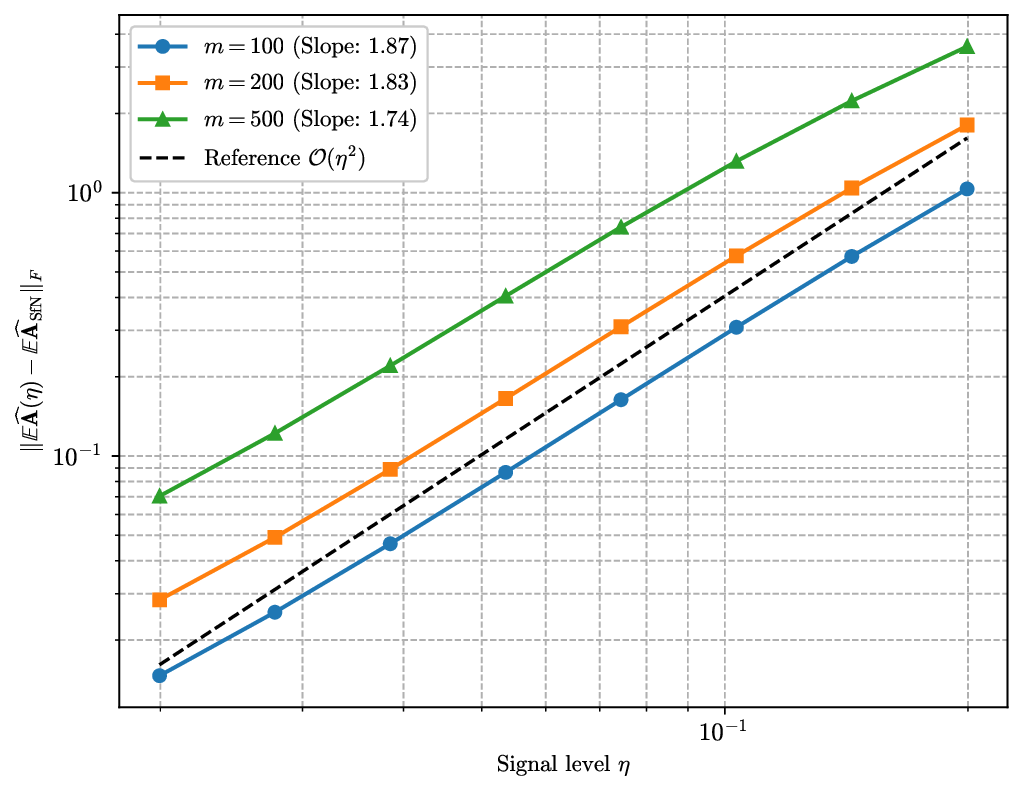}
\caption{Monte Carlos estimates of the population discrepancy $\|\E[\widehat{\BA}(\eta)]-\E \widehat{\BA}_{\rm SfN}\|_F$ versus signal strength $\eta$, on log--log axes, using $n=10^5$. The slope does not exactly equal 2 because of the slow convergence of the Monte Carlos simulation to estimate $\E\widehat{\BA}(\eta)$, which contributes to the first order term of $\eta$ and does not have a closed form expression.}
\label{fig:lsnr_expansion}
\end{figure}

\paragraph{Numerical illustration}
Figure~\ref{fig:lsnr_expansion} compares Monte Carlo estimates of the population discrepancy for Mickey Mouse-shaped templates and signals, using $m\in\{100,200,500\}$. The fitted log--log slopes of $\|\E\widehat{\BA}(\eta) - \E\widehat{\BA}_{\rm SfN}\|_F$ against $\eta$ are $1.87$, $1.83$, and $1.74$, respectively. For fixed dimensions, the same expansion order holds in the Frobenius norm. A nonzero quadratic leading term predicts an asymptotic log--log slope of two for fixed model parameters. Finite signal levels and Monte Carlo error can affect the fitted slopes, so the plotted values illustrate the trend without determining the asymptotic rate.

\paragraph{Conclusion}
We have characterized the structure induced by template matching and averaging in the generalized Procrustes model: the pure-noise population mean preserves the template's singular-vector bases while its spectrum may be distorted. The concentration bounds and weak-signal expansion quantify how this bias becomes visible as observations are averaged. It would be interesting to identify analogous invariants and bias mechanisms in other applications, especially the original cryo-EM models, where density projections and additional processing steps make the induced structure harder to characterize.

\section{Proofs}\label{s:proof}

 The proofs in this section focus on analyzing statistical properties of $\widehat{\BA}_{\rm SfN}$. In Section \ref{ss:proof1}, we characterize the model bias by justifying the singular value decomposition of the expectation $\E \widehat{\BA}_{\rm SfN}$ and derive lower and upper bounds for its singular values, establishing Theorem \ref{thm:main} and Corollary \ref{cor:svd_efn}. In addition, we derive upper bounds for the strength of the random deviation and finally combine the results to complete the proof of Theorem \ref{thm:error}. Section~\ref{pf_sec: LSNR} proves Theorem~\ref{thm:main_lsnr} in the low SNR regime.

\subsection{Pure-noise regime}\label{ss:proof1}


In this section, we assume that the observation $\BA_i = \BDelta_i$ is pure noise. We analyze the statistical properties of $\widehat{\BA}_{\rm SfN} = n^{-1}\sum_{i=1}^n \widehat{\BR}_i^{\top}\BA_i$. 

\paragraph{Gaussian coordinates and estimator decomposition}
Let $\BA^{\#}=\BU^{\#}\BSigma^{\#}\BV^{\#\top}$, where $\BU^{\#}\in\Od(d)$ and $\BV^{\#\top}\in\St(d,m)$. Define
\begin{equation}\label{eq: Gaussian_coordinates}
\BPhi_i:=\BDelta_i\BV^{\#}\in\mathbb{R}^{r\times d},
\qquad
\BDelta_{i,\perp}:=
\BDelta_i\big(\I_m-\BV^{\#}\BV^{\#\top}\big)
\end{equation}
and it holds that
\[
\BDelta_i = \BPhi_i \BV^{\#\top} + \BDelta_{i,\perp}.
\]
The matrices $\BPhi_i$ defined in~\eqref{eq: Gaussian_coordinates} are independent standard Gaussian matrices. For each $i$, the matrices $\BPhi_i$ and $\BDelta_{i,\perp}$ are jointly Gaussian and independent. 
Also, $\E[\BDelta_{i,\perp}]=\bzero$ holds. 

Define
\[
\BQ_i:=\widehat{\BR}_i\BU^{\#}\in\RR^{r\times d},~~~ \BQ_i \in\St(r,d)
\]
where $\widehat{\BR}_i\in\St(r,d)$, and they satisfy
\[
\begin{aligned}
\widehat{\BR}_i \in \argmax_{\BR \in \St(r,d)} \lag \BR\BA^{\#}, \BDelta_i \rag &  \Longleftrightarrow~~ \widehat{\BR}_i \in \argmax_{\BR \in \St(r,d)} \lag \BR\BU^{\#}\BSigma^{\#}, \BDelta_i\BV^{\#} \rag  \\
& \Longleftrightarrow~~  \BQ_i= \widehat{\BR}_i\BU^{\#}\in\argmax_{\BQ\in\St(r,d)}
\langle\BQ, \BPhi_i\BSigma^{\#} \rangle
\end{aligned}
\]
The matrix $\BQ_i$ depends on $\BPhi_i\BSigma^{\#}$. Here and below, $\PP(\BPhi_i\BSigma^{\#})$ denotes this induced maximizing rule in template coordinates; no equivariance of the original selection rule is required.
Using $\BDelta_i=\BPhi_i\BV^{\#\top}+\BDelta_{i,\perp}$ in~\eqref{eq: Gaussian_coordinates} and substituting into~\eqref{eq: estimator_EfN}, we obtain
\begin{equation}\label{eq:Aerror}
\begin{aligned}
\widehat{\BA}_{\rm SfN} & = \frac{1}{n} \sum_{i=1}^n \widehat{\BR}_i^{\top}\BDelta_i = \frac{1}{n} \sum_{i=1}^n \BU^{\#}\BQ^{\top}_i\BDelta_i   \\
& = \frac{1}{n} \sum_{i=1}^n \BU^{\#}\BQ^{\top}_i\underbrace{\BDelta_i\BV^{\#}}_{\BPhi_i}\BV^{\#\top} +  \frac{1}{n} \sum_{i=1}^n \BU^{\#}\BQ^{\top}_i\BDelta_i(\I_m - \BV^{\#}\BV^{\#\top} ) \\
& = \underbrace{\BU^{\#}\E[\BQ^{\top}_i\BPhi_i]\BV^{\#\top}}_{\text{population mean}} 
+\underbrace{\frac{1}{n}\sum_{i=1}^n
\BU^{\#}\big(\BQ^{\top}_i\BPhi_i-\E[\BQ^{\top}_i\BPhi_i]\big)\BV^{\#\top}}_
{\text{centered fluctuation in the template coordinates}} +\underbrace{\frac{1}{n}\sum_{i=1}^n
\BU^{\#}\BQ^{\top}_i\BDelta_{i,\perp}}_
{\text{orthogonal noise term}}.
\end{aligned}
\end{equation}
The second term is centered by construction. The third has mean zero because $\BQ_i$ and $\BDelta_{i,\perp}$ are independent and $\E\BDelta_{i,\perp}=\bzero$.

To prove Theorem~\ref{thm:main}, we analyze the population term in~\eqref{eq:Aerror} and use Gaussian sign symmetry to show that $\E[\BQ_i^{\top}\BPhi_i]$ is diagonal.

\begin{proof}[\bf Proof of Theorem~\ref{thm:main}]

Based on the decomposition~\eqref{eq:Aerror}, we note that
\[
\BSigma_{\rm SfN} : = \E [\BQ^{\top}  \BPhi]  = \E [ \PP(\BPhi\BSigma^{\#})^{\top} \BPhi ]
\]
where $\BPhi\in\RR^{r\times d}$ has independent standard Gaussian entries.

We assume the rank of $\BA^{\#}$ equals $k$ and then 
\[
\BSigma^{\#}
=\blkdiag\left(\BSigma^{\#}_+,\bzero_{(d-k)\times(d-k)}\right),
\qquad \BSigma^{\#}_+=\diag(\sigma_1,\cdots,\sigma_k),
\]
where $\sigma_1\geq\cdots\geq\sigma_k>0$. For a generic copy $\BPhi$ of the matrix $\BPhi_i$ in~\eqref{eq: Gaussian_coordinates}, we let $\BPhi_+\in\RR^{r\times k}$ be the first $k$ columns of $\BPhi$, and write $\BPhi_-$ for its remaining $d-k$ columns.

Since $\BQ=\PP(\BPhi\BSigma^{\#})\in\St(r,d)$, we have $\BQ\BQ^{\top}=\PP(\BPhi\BSigma^{\#})\PP(\BPhi\BSigma^{\#})^{\top}=\I_r$. 
If $r\leq k\leq d$, then $\rank(\BPhi\BSigma^{\#})=r$ almost surely, and the polar factor is unique. We have
\[
\BPhi = [\BPhi_+, \BPhi_-],~~~\underbrace{\BPhi\BSigma^{\#}}_{\RR^{r\times d}} = [\underbrace{\BPhi_+ \BSigma_+^{\#}}_{\RR^{r\times k}}, \bzero_{r\times (d-k)}].
\]
Then we have
\[
\PP(\BPhi\BSigma^{\#}) = [\PP(\BPhi_+ \BSigma_+^{\#}), 0] \in\St(r,d)
\]
for $d\geq k\geq r$. Therefore, we have
\[
\BSigma_{\rm SfN} = \E [\PP(\BPhi\BSigma^{\#})^{\top}\BPhi] = 
\E \begin{bmatrix}
\PP(\BPhi_+ \BSigma_+^{\#})^{\top}\BPhi_+ & \PP(\BPhi_+ \BSigma_+^{\#})^{\top}\BPhi_- \\
0 & 0 
\end{bmatrix} =  \begin{bmatrix}
\E[\PP(\BPhi_+ \BSigma_+^{\#})^{\top}\BPhi_+] & 0 \\
0 & 0 
\end{bmatrix} 
\]
Otherwise, if $k < r$, then
\[
\BQ = \PP(\BPhi\BSigma^{\#}) = \PP\left( [\BPhi_+ \BSigma_+^{\#}, \bzero_{r\times (d-k)}]\right) 
= [\underbrace{\PP(\BSigma_+^{\#}\BPhi_+^{\top})^{\top}}_{\RR^{r\times k}}, \underbrace{\BE}_{\RR^{r\times (d-k)}}] \in \St(r,d)
\]
where $\BPhi_+\BSigma_+^{\#}\in\RR^{r\times k}$ with $k\leq r$ and $\BE\in\RR^{r\times (d-k)}$, and they satisfy
\[
\PP(\BSigma^{\#}_+ \BPhi_+^{\top} )^{\top}\PP(\BSigma^{\#}_+ \BPhi_+^{\top} ) + \BE\BE^{\top} = \I_r,~~~\PP(\BSigma^{\#}_+ \BPhi_+^{\top} )\BE = \bzero_{k\times (d-k)},~~~\BPhi_+^{\top}\BE = 0.
\]

Then
\begin{align*}
\E[\BQ^{\top}\BPhi] & = \E\left[ \PP(\BPhi\BSigma^{\#})^{\top}\BPhi\right]  = \E\left[ \begin{bmatrix}
\PP(\BSigma_+^{\#}\BPhi_+^{\top}) \\
\BE^{\top}
\end{bmatrix} \BPhi\right] \\
& = \E \begin{bmatrix}
\PP(\BSigma_+^{\#}\BPhi_+^{\top}) \BPhi_+ & \PP(\BSigma_+^{\#}\BPhi_+^{\top}) \BPhi_- \\
\BE^{\top}\BPhi_+ & \BE^{\top}\BPhi_-
\end{bmatrix} = 
\begin{bmatrix}
\E \PP(\BSigma_+^{\#}\BPhi_+^{\top}) \BPhi_+ & \bzero_{k\times (d-k)} \\
\bzero_{(d-k)\times k} & \bzero_{(d-k)\times (d-k)} 
\end{bmatrix}
\end{align*}
The completion $\BE$ depends only on $\BPhi_+\BSigma^{\#}_+$ and is independent of $\BPhi_-$. Since $\BE^{\top}\BPhi_+=0$, we have $\E[\BE^{\top}\BPhi_-]=0$.

For both rank regimes, define
\[
\begin{aligned}
\BLambda_+ :=\E\left[\PP(\BSigma_+^{\#}\BPhi_+^{\top})\BPhi_+\BSigma_+^{\#}\right] = \E\left[\big(\BSigma_+^{\#}\BPhi_+^{\top}\BPhi_+\BSigma_+^{\#}\big)^{1/2}\right]
\in\RR^{k\times k}.
\end{aligned}
\]
The second equality follows from Proposition~\ref{prop:polar}. For every fixed $\bv\ne\bzero$, $\BPhi_+\BSigma_+^{\#}\bv\ne\bzero$ almost surely, so $\bv^{\top}\BLambda_+\bv>0$. Thus $\BLambda_+\succ0$, including when $k>r$.

For any diagonal sign matrix $\BD=\diag(\eps_1,\cdots,\eps_k)$, with $\eps_j\in\{-1,1\}$, Gaussian symmetry gives $\BPhi_+\BD\stackrel{d}{=}\BPhi_+$. We therefore obtain the following chain of equalities for the expectation:
\[
\begin{aligned}
\BLambda_+
&= \E \left[
\big(\BSigma_+^{\#}(\BPhi_+\BD)^{\top}(\BPhi_+\BD)\BSigma_+^{\#}\big)^{1/2}
\right]\\
&= \E \left[ \big(\BSigma_+^{\#}\BD\BPhi_+^{\top}\BPhi_+\BD\BSigma_+^{\#}\big)^{1/2}
\right]\\
&= \E \left[ \big(\BD\BSigma_+^{\#}\BPhi_+^{\top}\BPhi_+\BSigma_+^{\#}\BD\big)^{1/2}
\right]\\
&=\BD  \E \left[ \big(\BSigma_+^{\#}\BPhi_+^{\top}\BPhi_+\BSigma_+^{\#}\big)^{1/2}
\right]\BD \\
&=\BD\BLambda_+\BD.
\end{aligned}
\]
The first equality uses equality in distribution. The third uses the commutativity of $\BD$ and $\BSigma^{\#}_+$, and the fourth uses orthogonality of $\BD$ and uniqueness of the positive semidefinite square root. Thus, $\BLambda_+=\BD\BLambda_+\BD$ for every such $\BD$ because for $a\ne b$, choosing $\eps_a=-1$ and all other signs equal to $1$ gives $\Lambda_{ab}=-\Lambda_{ab}$, so $\Lambda_{ab}=0$.



Therefore, it holds that
\[
\E[\BQ^{\top}\BPhi]=\BSigma_{\rm SfN}=\blkdiag\big(\BLambda_+\BSigma^{\#-1}_+,
\bzero_{(d-k)\times(d-k)}\big)
\]
is diagonal, with positive first $k$ entries, for every $1\leq k\leq d$.
\end{proof}

Taking expectations in~\eqref{eq:Aerror} now yields~\eqref{eq: population_mean_EfN}. Since $\BU^{\#}$ and $\BV^{\#}$ have orthonormal columns and $\BSigma_{\rm SfN}$ is diagonal and nonnegative, this is a singular-value decomposition of the population mean with the same left and right singular-vector bases as the template.

\begin{proof}[\bf Proof of Corollary~\ref{cor:svd_efn}]

By Theorem~\ref{thm:main}, it holds that
\[
\BLambda_+
= \E \left[
\big(\BSigma^{\#}_+\BPhi_+^{\top}\BPhi_+\BSigma^{\#}_+\big)^{1/2}
\right],
\]
where $\BPhi_+\in\mathbb{R}^{r\times k}$ is standard Gaussian. In Step~1, we use orthogonal invariance to compute $\E[(\BPhi_+^{\top}\BPhi_+)^{1/2}]$. In Step~2, we use this expectation and the extreme template singular values $\sigma_1$ and $\sigma_k$ to derive lower and upper bounds for $\BLambda_+$.

\noindent{\bf Step 1: Orthogonal invariance of the Gaussian square root.}
For every $\BO\in\Od(k)$, we have $\BPhi_+\BO\stackrel{d}{=}\BPhi_+$. Hence,
\[
\begin{aligned}
 \E \left[(\BPhi_+^{\top}\BPhi_+)^{1/2}\right]
&= \E \left[
\big((\BPhi_+\BO)^{\top}(\BPhi_+\BO)\big)^{1/2}
\right] =\BO^{\top} 
 \E \left[(\BPhi_+^{\top}\BPhi_+)^{1/2}\right]\BO.
\end{aligned}
\]
Taking $\BO$ to be a diagonal sign matrix shows that this expectation is diagonal; taking $\BO$ to be a permutation matrix shows that all its diagonal entries are equal. Taking traces identifies the common value as $k^{-1} \E \|\BPhi_+\|_*$. Thus, using Proposition~\ref{prop:polar},
\begin{equation}\label{eq: def_c_rd}
 \E \left[\PP(\BPhi_+^{\top})\BPhi_+\right]
= \E \left[(\BPhi_+^{\top}\BPhi_+)^{1/2}\right]
=\frac{ \E \|\BPhi_+\|_*}{k}\I_k.
\end{equation}

\noindent{\bf Step 2: Comparison with the template spectrum.}
For every rectangular matrix $\BX$, the SVD gives the identity
\[
(\BX^{\top}\BX)^{1/2}=\BX^{\top}\bigl[(\BX\BX^{\top})^{\dagger}\bigr]^{1/2}\BX.
\]
We apply this identity with $\BX=\BPhi_+\BSigma_+^{\#}$.
Since $\sigma_k^2\I_k\preceq\BSigma^{\# 2}_+\preceq\sigma_1^2\I_k$, congruence by $\BPhi_+$ yields
\[
\sigma_k^2\BPhi_+\BPhi_+^{\top}
\preceq\BPhi_+\BSigma^{\# 2}_+\BPhi_+^{\top}
\preceq\sigma_1^2\BPhi_+\BPhi_+^{\top}.
\]
The matrices $\BPhi_+\BPhi_+^{\top}$ and $\BPhi_+\BSigma^{\# 2}_+\BPhi_+^{\top}$ have the same range $\Ran(\BPhi_+)$ and they are positive definite on this common range. Therefore, pseudo-inverse reverses the order there, and the positive semidefinite square root preserves it: 
\[
\frac{1}{\sigma_1}
\bigl[(\BPhi_+\BPhi_+^{\top})^{\dagger}\bigr]^{1/2}
\preceq
\bigl[(\BPhi_+\BSigma^{\# 2}_+\BPhi_+^{\top})^{\dagger}\bigr]^{1/2}
\preceq
\frac{1}{\sigma_k}
\bigl[(\BPhi_+\BPhi_+^{\top})^{\dagger}\bigr]^{1/2}.
\]
Multiplying the inequality on the left by $\BSigma^{\#}_+\BPhi_+^{\top}$ and on the right by $\BPhi_+\BSigma^{\#}_+$, then using the identity gives
\[
\frac{1}{\sigma_1}\BSigma^{\#}_+
(\BPhi_+^{\top}\BPhi_+)^{1/2}\BSigma^{\#}_+
\preceq
\big(\BSigma^{\#}_+\BPhi_+^{\top}\BPhi_+\BSigma^{\#}_+\big)^{1/2}
\preceq
\frac{1}{\sigma_k}\BSigma^{\#}_+
(\BPhi_+^{\top}\BPhi_+)^{1/2}\BSigma^{\#}_+.
\]
Taking expectations and using~\eqref{eq: def_c_rd} gives
\begin{equation}\label{eq:Lambda_+}
\frac{ \E \|\BPhi_+\|_*}{k\sigma_1}\BSigma^{\# 2}_+
\preceq\BLambda_+
\preceq\frac{ \E \|\BPhi_+\|_*}{k\sigma_k}\BSigma^{\# 2}_+.
\end{equation}
Combining this with $\sigma_k^2\I_k\preceq\BSigma^{\# 2}_+\preceq\sigma_1^2\I_k$ proves~\eqref{eq: singular_value_bounds_Lambda}. Since $\BLambda_+$ and $\BSigma^{\#}_+$ are diagonal, dividing their $j$th diagonal entries by $\sigma_j>0$ gives~\eqref{eq: singular_value_bounds_tau}. Finally, when $\BSigma^{\#}_+=\I_k$, then  
\[
\BLambda_+= \frac{ \E \|\BPhi_+\|_*}{k}\I_k,
\]
proving~\eqref{eq: singular_value_isotropic_EfN}.

\noindent{\bf Step 3: Lower bound of $\| \E \widehat{\BA}_{\rm SfN}\|$.}
We next derive a lower bound for the operator norm of the population mean, the first term in~\eqref{eq:Aerror}. 
%
It is obvious that
\[
\begin{aligned}
 \frac{\E\|\BPhi_+\|_*}{k\sigma_1}  \BSigma_+^{\#}\preceq \BLambda_+ \BSigma_+^{\# -1} \preceq \frac{\E\|\BPhi_+\|_*}{k\sigma_k}  \BSigma_+^{\#},~~~ \left\|\E[\widehat{\BA}_{\rm SfN}]\right\|
 =\left\|\BLambda_+\BSigma^{\#-1}_+\right\| \geq \frac{\E\|\BPhi_+\|_*}{k} 
\end{aligned}
\]
which follows from~\eqref{eq:Lambda_+}.

%

For the Gaussian nuclear norm, the singular values satisfy
\[
\begin{aligned}
\|\BPhi_+\|_F^2 &  =  \sum_{j=1}^{r\wedge k} \sigma_j^{2/3}(\BPhi_+)\cdot \sigma_j^{4/3} (\BPhi_+) \leq \left(\sum_{j=1}^{r\wedge k} \sigma_j(\BPhi_+)\right)^{2/3} \left(\sum_{j=1}^{r\wedge k} \sigma_j^4(\BPhi_+)\right)^{1/3}
\end{aligned}
\]
which follows from the H\"older inequality. Using the H\"older inequality again after taking the expectation, we have
\[
\E \|\BPhi_+\|_F^2 \leq \left[\E \left(\sum_{j=1}^{r\wedge k} \sigma_j(\BPhi_+)\right) \right]^{2/3} \cdot
\left[  \E \sum_{j=1}^{r\wedge k} \sigma_j^4(\BPhi_+) \right]^{1/3} = [\E \|\BPhi_+\|_*]^{2/3} \cdot [\E \Tr\big[(\BPhi_+^{\top}\BPhi_+)^2\big]]^{1/3}.
 \]
Here 
\begin{align*}
\E \|\BPhi_+\|_F^2  = rk, \quad \E \Tr\big[(\BPhi_+^{\top}\BPhi_+)^2\big] = k(r(r+2)) + k(k-1)r = rk (r+k+1)
\end{align*}
Here each diagonal entry of $\BPhi_+^{\top}\BPhi_+$ has a $\chi_r^2$ distribution, with second moment $r(r+2)$, and each off-diagonal entry has second moment $r$.
Therefore, we have
\[
\begin{aligned}
 \E \|\BPhi_+\|_* & \geq \frac{[\E\|\BPhi_+\|_F^2]^{3/2}}{[\E \Tr\big[(\BPhi_+^{\top}\BPhi_+)^2\big]]^{1/2} } = \frac{(rk)^{3/2}}{ \sqrt{rk(r+k+1)}  } = \frac{rk}{\sqrt{r+k+1}}, \\ 
 \left\|\E[\widehat{\BA}_{\rm SfN}]\right\| & \geq  \frac{r}{\sqrt{r+k+1}}.
 \end{aligned}
\]
%
\end{proof}


Next, we bound the second term in~\eqref{eq:Aerror} using the matrix Bernstein inequality with truncation.

\begin{lemma}[\bf Estimation of the second term in~\eqref{eq:Aerror}]\label{lemma: ub_signal_deviation}
Suppose $n\ge d\vee2$, $\gamma\ge1$, and $(\gamma+1)\log n\le\sqrt{nd}$. For the second term in~\eqref{eq:Aerror}, with probability at least $1-O(n^{-\gamma})$,
\begin{equation}\label{eq: signal_deviation_simplified}
\left\|\frac{1}{n}\sum_{i=1}^n
\BU^{\#}\big(\BQ^{\top}_i\BPhi_i-\E[\BQ^{\top}_i\BPhi_i]\big) \BV^{\#\top}\right\| \leq  C(\sqrt{r} + \sqrt{d})\sqrt{\frac{2\gamma\log n}{n}}
\end{equation}
where $\BQ_i=\PP(\BPhi_i\BSigma^{\#})$.
\end{lemma}
Before proving the lemma above, we state the matrix Bernstein inequality used in the proof.

\begin{theorem}[Matrix Bernstein inequality, {\cite[Theorem~6.1.1]{T15}}]\label{thm:bernstein}
Let $\{\BZ_i\}$ be a finite sequence of independent random matrices in $\mathbb{R}^{d_1\times d_2}$ satisfying
\[
\E[\BZ_i]=\bzero, \quad \|\BZ_i\|\leq R \quad\text{almost surely for every }i.
\]
Define
\[
v:=\max\left\{
\left\|\sum_i\E[\BZ_i\BZ_i^{\top}]\right\|,
\left\|\sum_i\E[\BZ_i^{\top}\BZ_i]\right\|
\right\}.
\]
Then, for every $t>0$,
\[
\Pr \left\{\left\|\sum_i\BZ_i\right\|\geq t\right\}
\le(d_1+d_2)\exp \left(-\frac{t^2/2}{v+Rt/3}\right).
\]
\end{theorem}

\begin{proof}[\bf Proof of Lemma~\ref{lemma: ub_signal_deviation}]
Define
\[
\BZ_i:=\BQ^{\top}_i\BPhi_i-\E[\BQ^{\top}_i\BPhi_i]
\in\mathbb{R}^{d\times d}.
\]
The $\BZ_i$ are independent, identically distributed, and centered.  We first control the summand norms and the matrix variance, then apply Theorem~\ref{thm:bernstein} to centered, truncated summands.

\noindent\textbf{Step 1: Summand norms and matrix variance.}
Every fitted factor $\BQ_i = \PP(\BPhi_i\BSigma^{\#})$ belongs to $\St(r,d)$, so $\BQ_i\BQ^{\top}_i=\I_r$. 

 By the triangle inequality, and the Gaussian norm estimate in~\cite[Theorem~7.3.1]{V19},
\[
\begin{aligned}
\|\BZ_i\| &\leq\|\BPhi_i\| +\bigl\|\E[\BQ^{\top}_i\BPhi_i]\bigr\| \leq\|\BPhi_i\| + \E \|\BPhi_i\| \leq\|\BPhi_i\|+\sqrt{r}+\sqrt{d}  \leq 2(\sqrt{r} + \sqrt{d}) + t
\end{aligned}
\]
with probability at least $1-2e^{-t^2/2}$, since Gaussian concentration~\cite[Proposition~10.3]{HMT11} yields
\begin{equation}\label{eq: signal_deviation_norm_tail}
\Pr \left\{\|\BPhi_i\|>\sqrt{r}+\sqrt{d}+t\right\}
\leq 2e^{-t^2/2},\qquad t\geq 0.
\end{equation}

For the variance bound, the same isometry identity yields
\[
\begin{aligned}
\E[\BZ_i^{\top}\BZ_i]
& = \E[\BPhi_i^{\top}\BPhi_i] -\E[\BQ^{\top}_i\BPhi_i]^{\top}\E[\BQ^{\top}_i\BPhi_i] =r\I_d- \E[\BQ^{\top}_i\BPhi_i]^{\top}\E[\BQ^{\top}_i\BPhi_i]
\preceq r\I_d, \\
\E[\BZ_i\BZ_i^{\top}] & = \E[\BQ^{\top}_i\BPhi_i\BPhi_i^{\top}\BQ_i]-\E[\BQ^{\top}_i\BPhi_i]\left[\E[\BQ^{\top}_i\BPhi_i]\right]^{\top} \preceq \E \|\BPhi_i\|^2 \I_d \preceq C(r+d)\I_d.
\end{aligned}
\]
The last inequality follows by integrating the Gaussian norm tail. Hence
\begin{equation}\label{eq: second_moment_condition}
\max\left\{
\left\|\sum_{i=1}^n\E[\BZ_i\BZ_i^{\top}]\right\|,
\left\|\sum_{i=1}^n\E[\BZ_i^{\top}\BZ_i]\right\|
\right\}\leq Cn(r+d).
\end{equation}

\paragraph{Step 2: Truncation and the matrix Bernstein inequality.}
Let
\[
R=2(\sqrt r+\sqrt d)+\sqrt{2(\gamma+1)\log n},
\qquad
\widetilde{\BZ}_i
=\BZ_i1_{\{\|\BZ_i\|\leq R\}}
-\E[\BZ_i1_{\{\|\BZ_i\|\leq R\}}].
\]
The $\widetilde{\BZ}_i$ are independent and centered, have norm at most $2R$, and satisfy the same variance bound as the $\BZ_i$. By~\eqref{eq: signal_deviation_norm_tail}, a union bound and tail integration give
\[
\begin{aligned}
\Pr\{\max_i\|\BZ_i\|>R\} & \leq 2n^{-\gamma}, \\
\left\|\E[\BZ_i1_{\{\|\BZ_i\|\leq R\}}]\right\|
& = \left\|\E[\BZ_i1_{\{\|\BZ_i\|>R\}}]\right\| 
\leq \E\left[\|\BZ_i\|1_{\{\|\BZ_i\|>R\}}\right] \\
&\leq R\Pr(\|\BZ_i\|>R)+\int_R^{\infty}\Pr(\|\BZ_i\|>t)\diff t
\leq CRn^{-(\gamma+1)}.
\end{aligned}
\]
The bounds follow from centering and the Gaussian norm tail. Applying the matrix Bernstein inequality to $\widetilde{\BZ}_i$ and accounting for the truncation bias bounds $n^{-1}\|\sum_i\BZ_i\|$ by
\[
C\left( \sqrt{\frac{(r+d)(\log(2d)+\gamma\log n)}{n}}
+\frac{R(\log(2d)+\gamma\log n)}{n} +Rn^{-(\gamma+1)} \right)
\]
with probability at least $1-O(n^{-\gamma})$. Under the stated assumptions, $\log(2d)\lesssim\gamma\log n$ and
\[
\frac{R\sqrt{\gamma\log n}}{\sqrt{nd}}
\lesssim\sqrt{\frac{\gamma\log n}{n}}
+\frac{(\gamma+1)\log n}{\sqrt{nd}}\lesssim1.
\]
The last two terms are therefore absorbed into the first, proving~\eqref{eq: signal_deviation_simplified}.
\end{proof}

\begin{lemma}[\bf Estimation of the third term in~\eqref{eq:Aerror}]\label{lemma: ub_noise}
Under the assumptions of Theorem~\ref{thm:main}, for any $\gamma>0$, the following holds with probability at least $1-n^{-\gamma}$:
\begin{equation}\label{eq: orthogonal_noise_bound}
\left\|\frac{1}{n}\sum_{i=1}^n
\BU^{\#}\BQ^{\top}_i\BDelta_{i,\perp}\right\|
\leq \frac{\sqrt{r}+\sqrt{m-d}+\sqrt{2\gamma\log n}}{\sqrt n}.
\end{equation}
If $m=d$, the orthogonal noise term is identically zero.
\end{lemma}

\begin{proof}[\bf Proof of Lemma~\ref{lemma: ub_noise}]
If $m=d$, the orthogonal noise term is zero, so suppose $m-d\geq1$. By the Gaussianity of $\BDelta_i$, the matrices $\BPhi_i= \BDelta_i\BV^{\#}$ and $\BDelta_i(\I_m-\BV^{\#}(\BV^{\#})^{\top})$ are independent. Since the fitted factor $\BQ^{\top}_i$ depends only on $\BPhi_i=\BDelta_i\BV^{\#}\in\RR^{r\times d}$ under the prescribed selection rule, it is independent of $\BDelta_{i,\perp}$.

We choose $\BV_{\perp}^{\#}\in\mathbb{R}^{m\times (m-d)}$ so that $\BV_{\perp}^{\#}\BV_{\perp}^{\#\top} = \I_m - \BV^{\#}\BV^{\#\top}$ and then 
\[
\frac{1}{n}\sum_{i=1}^n\BU^{\#}\BQ^{\top}_i\BDelta_{i,\perp}  = \frac{1}{n} \BU^{\#}\sum_{i=1}^n\BQ^{\top}_i \BDelta_i (\I_m - \BV^{\#}\BV^{\#\top}) =  \frac{1}{n} \BU^{\#}\left[\sum_{i=1}^n\BQ^{\top}_i \BDelta_i \BV_{\perp}^{\#}\right]\BV_{\perp}^{\#\top}
\]
where $\BQ^{\top}_i = \PP( \BSigma^{\#} \BPhi_i^{\top})\in\RR^{d\times r}$ and $\BDelta_{i,\perp} = \BDelta_i (\I_m - \BV^{\#}\BV^{\#\top}).$
Therefore, it suffices to bound the operator norm of
\[
[\BQ^{\top}_1,\cdots, \BQ^{\top}_n] 
\begin{bmatrix}
\BDelta_1\BV^{\#}_{\perp} \\
\vdots \\
\BDelta_n\BV^{\#}_{\perp} 
\end{bmatrix} = \sum_{i=1}^n \BQ^{\top}_i \BDelta_i \BV_{\perp}^{\#}.
\]
Conditionally on $(\BPhi_1,\cdots,\BPhi_n)$, the matrix $\{\BQ_i\}_{i=1}^n$ are fixed, while the stacked residual matrix is standard Gaussian in $\RR^{nr\times (m-d)}$ and independent of $\BQ_1,\cdots,\BQ_n$. Since
\[
\|[\BQ_1^{\top},\cdots,\BQ_n^{\top}]\|\le\sqrt n,
\qquad
\|[\BQ_1^{\top},\cdots,\BQ_n^{\top}]\|_F=\sqrt{nr},
\]
Proposition~10.1 in~\cite{HMT11} gives
\[
\E\left[ \left\|\sum_i\BQ_i^{\top}\BDelta_i\BV_{\perp}^{\#}\right\| \middle|\,\BPhi_1,\cdots,\BPhi_n \right]\le\sqrt{n(m-d)}+\sqrt{nr}.
\]
Note that
$n^{-1}\|\sum_i\BQ_i^{\top}\BDelta_i\BV_{\perp}^{\#}\|$ is $n^{-1/2}$-Lipschitz with respect to the stacked Gaussian matrix in the Frobenius norm. The Gaussian concentration inequality~\cite[Proposition~10.3]{HMT11} therefore gives
\[
\Pr\left\{ \frac{1}{n}\left\|\sum_i\BQ_i^{\top}\BDelta_i\BV_{\perp}^{\#}\right\|
>\frac{\sqrt r+\sqrt{m-d}+t}{\sqrt n}
\Big| \BPhi_1,\cdots,\BPhi_n
\right\}\leq e^{-t^2/2}.
\]
Averaging over the conditioning variables and setting $t=\sqrt{2\gamma\log n}$ completes the proof.
\end{proof}


\begin{proof}[\bf Proof of Theorem~\ref{thm:error}]
Corollary~\ref{cor:svd_efn} gives the deterministic bound
\[
\left\|\E[\widehat{\BA}_{\rm SfN}]\right\|
\geq\frac{r}{\sqrt{r+k+1}}.
\]
Combining Lemmas~\ref{lemma: ub_signal_deviation} and~\ref{lemma: ub_noise} in~\eqref{eq:Aerror} gives~\eqref{eq: order_EfN_estimator} with probability at least $1-O(n^{-\gamma})$, after absorbing $\sqrt{2\gamma\log n}$ into $C(\sqrt r+\sqrt d)\sqrt{2\gamma\log n}$. The stated sample-size conditions follow from the displayed lower bound.
\end{proof}

\subsection{Low-SNR regime}\label{pf_sec: LSNR}
\begin{proof}[\bf Proof of Theorem~\ref{thm:main_lsnr}]
We use the notation of Theorem~\ref{thm:main_lsnr}.
\paragraph{Step 1: Expansion of the population mean.}
Recall that ${\cal T}(\BY)=\PP(\BY\BA^{\#\top})^{\top}\BY$. Then
\[
\E\widehat{\BA}(\eta)  = \frac{1}{n}\sum_{i=1}^n \E {\cal T}(\eta \BR_i^{\star}\BA^{\star} + \BDelta_i ).
\]
For the standard Gaussian density $\varphi(\cdot)$ on $\mathbb{R}^{r\times m}$, it holds that
\[
\varphi(\BY-\eta\BR_i^{\star}\BA^{\star})
=\varphi(\BY) \exp \left(\eta\langle\BY,\BR^{\star}_i\BA^{\star}\rangle -\frac{\eta^2}{2}\|\BR^{\star}_i\BA^{\star}\|_F^2\right).
\]
Changing variables in the Gaussian integral yields the exact identity
\begin{equation}\label{eq: LSNR_Gaussian_shift}
\E\widehat{\BA}(\eta) 
= \frac{1}{n}\sum_{i=1}^n \E \left[ {\cal T}(\BDelta_i)
\exp \left(\eta\langle\BDelta_i,\BR_i^{\star}\BA^{\star}\rangle-\frac{\eta^2}{2}\|\BR_i^{\star}\BA^{\star}\|_F^2\right) \right].
\end{equation}

Note that
\begin{align*}
\E\widehat{\BA}(0) & = \E \widehat{\BA}_{\rm SfN}, \\
\frac{\diff}{\diff \eta}\E \widehat{\BA}(\eta) \Big|_{\eta =0} & = \frac{1}{n} \sum_{i=1}^n \E {\cal T}(\BDelta_i) \lag \BDelta_i, \BR_i^{\star}\BA^{\star}\rag, \\
\frac{\diff^2}{\diff \eta^2}\E\widehat{\BA}(\eta) \Big|_{\eta =0} & = \frac{1}{n}\sum_{i=1}^n \E {\cal T}(\BDelta_i) (\lag \BDelta_i, \BR_i^{\star}\BA^{\star}\rag^2 - \|\BR_i^{\star}\BA^{\star}\|_F^2  ),\\
\frac{\diff^3}{\diff\eta^3}\E\widehat{\BA}(\eta)\Big|_{\eta=0} & = \frac{1}{n}\sum_{i=1}^n\E
{\cal T}(\BDelta_i)\big(\langle\BDelta_i,\BR_i^{\star}\BA^{\star}\rangle^3 - 3\|\BR_i^{\star}\BA^{\star}\|_F^2 \langle\BDelta_i,\BR_i^{\star}\BA^{\star}\rangle\big).
\end{align*}

Since $\|{\cal T}(\BY)\|_F=\|\BY\|_F$, the differentiated integrands through order four are dominated on compact $\eta$-intervals by a Gaussian-integrable function of the form \[ C\|\BY\|_F(1+\|\BY\|_F^4)e^{C\|\BY\|_F}\varphi(\BY). \] Thus the differentiation under the integral and Taylor's theorem apply.

\paragraph{Step 2: Sign symmetry and the quadratic correction.}
Suppose $k\geq r$. For a full row-rank matrix $\BX$, its polar equals
\begin{equation}\label{eq: odd_polar_retraction}
\PP(\BX)=(\BX\BX^{\top})^{-1/2}\BX,
\qquad
\PP(-\BX)=-\PP(\BX).
\end{equation}
When $k\ge r$, $\BY\BA^{\#\top}$ has full row rank almost surely for a generic Gaussian matrix $\BY$, and hence $ {\cal T}(-\BY)={\cal T}(\BY)$ holds. In addition, we note that the Gaussian symmetry implies $\E\widehat{\BA}(-\eta)=\E\widehat{\BA}(\eta)$. Thus its first and third derivatives at zero vanish, and the bounded fourth derivative gives~\eqref{eq: LSNR_quadratic_expansion}.
%
%
%
\end{proof}

\section*{AI Disclaimer}
The authors acknowledge the use of \textbf{Gemini 3.1 Pro, Gemini 3.6 Thinking} (developed by \textbf{Google}) and \textbf{ChatGPT 6 Sol, ChatGPT 6 Astra and ChatGPT 5.6 Sol} (developed by \textbf{OpenAI}) for language polishing, proofreading, literature searches, checking mathematical arguments, and improving figure presentation during the manuscript preparation process. All mathematical derivations, scientific content, and final text were independently verified, calculated, and approved by the human authors, who hold sole responsibility for the work.

\bibliography{reference.bib}
\bibliographystyle{abbrv}

\end{document}